\documentstyle[12pt]{article}
\title{Hopf Algebras with Positive Bases}
\begin{document}
\author{Jiang-Hua Lu\thanks{Research is supported by . }\\
Department of Mathematics\\ 
University of Arizona\\
\\
Yong-Chang Zhu\thanks{Research is supported by . },
\hspace{5mm} Min Yan \\
Department of Mathematics\\ 
Hong Kong University of Science and Technology\\
Clear Water Bay, Hong Kong} 
\date{}
\maketitle

\newcommand{\la}{\leftarrow}
\newcommand{\ra}{\rightarrow}
\newcommand{\da}{\downarrow}
\newcommand{\ua}{\uparrow}
\newcommand{\uda}{\updownarrow}
\newcommand{\sub}{\subset}
\newcommand{\lra}{\longrightarrow}
\newcommand{\lla}{\longleftarrow}
\newcommand{\La}{\Leftarrow}
\newcommand{\Ra}{\Rightarrow}
\newcommand{\Ua}{\Uparrow}
\newcommand{\Da}{\Downarrow}
\newcommand{\Uda}{\Dpdownarrow}
\newcommand{\Lla}{\Longleftarrow}
\newcommand{\Lra}{\Longrightarrow}
\newcommand{\Llra}{\Longleftrightarrow}
\newcommand{\pa}{\partial}
\newcommand{\epi}{\mbox{epi}}

\newtheorem{th}{Theorem}
\newtheorem{lem}{Lemma}
\newtheorem{cor}[lem]{Corollary}
\newtheorem{prop}[lem]{Proposition}
\newtheorem{df}{Definition}
\newtheorem{rmk}{Remark}
\newtheorem{eg}{Example}

\begin{abstract}
We show that if a finite dimensional Hopf algebra over ${\bf C}$ has a basis such
that all the structure constants are non-negative, then the Hopf algebra must be
given by a finite group $G$ and a factorization $G=G_+G_-$ into two subgroups. We also
show that Hopf algebras in the category of finite sets with correspondences as
morphisms are classified in the similar way. Our results can be used to explain some
results in Hopf algebras from set-theoretical viewpoint.
\end{abstract}

\tableofcontents

\section{Introduction}

   Recently much progress has been made in the classification of
  finite dimensional Hopf algebras (see \cite{EG} \cite{Ma} \cite{BDG} and 
  references quoted there). The classification
  theorems are often established under the assumption that
  the dimensions of the Hopf algebras have few divisors. 
  
  In this work, 
  we prove two classification theorems under different
  kinds of assumptions.
   The first theorem (Theorem \ref{main}) classifies 
  finite dimensional Hopf algebras over $\bf C$ 
  that admit positive bases.  
  A basis of a Hopf algebra over $\bf C$ is called positive if all the 
  structure coefficients under the basis are non-negative real numbers
  (see Section 2 for more precise definitions). This class of Hopf
  algebras is closed under tensor product, dual and the Drinfel'd
  double construction. Our theorem says that if a Hopf algebra has a positive basis,
then it must be isomorphic to the bicrossproduct Hopf algebra induced from a finite
group $G$ with a unique factorization $G=G_+G_-$ of two subgroups. The construction of
this type of Hopf algebras can be found in \cite{T} and \cite{mj:paper}.

 The second theorem (Theorem \ref{correspondence}) classifies Hopf algebras in the
correspondence category.
  In the monoidal category of finite sets, we may formally introduce the concept of
Hopf algebra by imposing all the axioms in terms of commutative diagrams. Such a
definition is just that of finite groups. 
   The correspondence category under our consideration still has finite sets as 
   objects. But a morphism
   from $X$ to $Y$ is a subset of $ X\times Y$ instead
   of a map from $X$ to $Y$. This is a monoidal category in which the concept of Hopf
algebra can be  similarly defined.
Our theorem shows that a Hopf algebra in the correspondence category is still given by
a finite group $G$ with a unique factorization $G=G_+G_-$ of two subgroups.

    One of the motivations to consider positive based
    Hopf algebras is that there are many natural Hopf algebras  with 
    nearly positive bases. We call a basis nearly positive if all the structure
coefficients
     except the antipode under the basis are non-negative
     real numbers. For example, the modified
     quantum group $U_q ( g )$ in \cite{Lusztig} has a nearly positive basis realized
in terms certain varieties. Another example is the direct sum 
$ \oplus_{1}^{\infty} R (S_n)$ and $ \oplus_{1}^{\infty} R (GL_n ( F_q ))$ of complex
representation rings \cite{Zelevinsky}, where positive bases are given by irreducible
representations.  The positivity of structural 
     coefficients are due to the very ways of defining the
   multiplication and the comultiplication. 
See \cite{K} for more 
    examples.  All the examples above are infinite dimensional
   Hopf algebras.   In the finite dimensional case, we know no example of a nearly
positive
   basis that is not a positive basis.

In a subsequent paper \cite{LYZ}, we classify positive quasi-triangular structures on
positively based Hopf algebras, and discuss its relation to set-theoretical
Yang-Baxter equations.

The paper is organized as follows: We first present the basic concepts and notations,
and state the first classification theorem. In the next two sections, we give detailed
proof of the theorem.  In the final section, we discuss the correspondence category
and outline the proof of the second classification theorem.

\section{Positively Based Hopf Algebras}

Given a group $G$, the group algebra ${\bf C}G$ and the dual group algebra 
$({\bf C}G)^*$ have $G$ and $G^*$ as positive bases in the sense that all the
structure constants are nonnegative. More generally, we call a basis $B$ of a
Hopf algebra $H$ {\em positive} if 
\begin{enumerate}
\item The coordinates of the unit $1$ are non-negative;
\item The coordinates of the counit $\epsilon$ (with respect to the dual basis $B^*$)
are non-negative, i.e., $\epsilon(b)\geq 0$ for all $b\in B$;
\item For any $b_1,b_2\in B$, the coordinates of $b_1b_2$ are non-negative;
\item For any $b\in B$, the coordinates of $\Delta(b)$ (with respect to the tensor
basis $B\otimes B$) are non-negative;
\item For any $b\in B$, the coordinates of the antipode $S(b)$ are non-negative.
\end{enumerate}
We call a Hopf algebra {\em positively based} if it has a positive basis.

It is easy to see that the dual and the tensor of positively
based Hopf algebras are still positively based. For the antipode $S$ of a finite
dimensional Hopf algebra, it is well-known that $S^N=1$ for some positive integer
$N$ \cite{R}. This implies that if $B$ is a positive basis of $H$, then the
coordinates of $S^{-1}=S^{N-1}$ with respect to $B$ are still
non-negative. Therefore $B^*\otimes B$ is a positive basis of the Drinfel'd double
$D(H)$.

The prototypical example of positively
based Hopf algebras is the bicrossproduct Hopf algebra $H(G;G_+,G_-)$ induced from a
finite group $G$ with a unique factorization $G=G_+G_-$ of two subgroups (see \cite{T}
and \cite{mj:paper}). By $G=G_+G_-$ we mean $G_{\pm}$ are subgroups of $G$, such that
any
$g\in G$ can be written as $g=g_+ g_-$ for unique $g_+\in G_+$ and $g_-\in G_-$. By
taking inverse, any $g\in G$ can also be written as
$g=\bar{g}_-\bar{g}_+$ for unique $\bar{g}_+\in G_+$ and
$\bar{g}_-\in G_-$. Therefore to any $g\in G$ is associated four uniquely
determined elements 
\[
\alpha_+(g)=g_+,\hspace{3mm} 
\beta_+(g)=\bar{g}_+, \hspace{3mm}
\alpha_-(g)=\bar{g}_-, \hspace{3mm}
\beta_-(g)=g_-, 
\]
given by
\begin{equation}\label{gpd5}
g=g_+ g_-=\bar{g}_-\bar{g}_+.
\end{equation}
It is then easy to see that
\begin{equation}\label{gpd6}
G_-\times G_+\ra G_+,\hspace{5mm}
(\bar{g}_-,\bar{g}_+) \mapsto (\bar{g}_-\bar{g}_+)_+=g_+,
\end{equation}
\begin{equation}\label{gpd7}
G_-\times G_+\ra G_-,\hspace{5mm}
(\bar{g}_-,\bar{g}_+) \mapsto (\bar{g}_-\bar{g}_+)_-=g_-.
\end{equation}
\begin{equation}\label{gpd8}
G_+\times G_-\ra G_+,\hspace{5mm}
(g_+,g_-) \mapsto \overline{(g_+g_-)}_+=\bar{g}_+,
\end{equation}
\begin{equation}\label{gpd9}
G_+\times G_-\ra G_-,\hspace{5mm}
(g_+,g_-) \mapsto \overline{(g_+g_-)}_-=\bar{g}_-
\end{equation}
are all group actions. Because of this, we also write
\[
g_+=\;^{\bar{g}_-}\bar{g}_+,\hspace{3mm} 
g_-=\bar{g}_-^{\;\bar{g}_+}, \hspace{3mm}
\bar{g}_+=g_+^{\;g_-}, \hspace{3mm}
\bar{g}_-=\;^{g_+}g_-.
\]
The third and the fourth actions satisfy the following matching relations
\begin{equation}\label{gpdmatch}
\;^{g_+}(g_-h_-)=\;^{g_+}g_-\;\;^{(g_+^{\;g_-})}h_-,\hspace{5mm}
(h_+g_+)^{g_-}=h_+^{(\;^{g_+}g_-)}\;g_+^{\;g_-}.
\end{equation}
It is well-known that given two actions between two groups $G_{\pm}$ satisfying the
matching relations, we may reconstruct $G$. Similarly, $G$ can be reconstructed from
the first and the second actions, which satisfy similar matching relations. 
We will not use the first and the second actions in the proof of
the classification theorem.

The Hopf algebra $H(G;G_+,G_-)$ has basis $\{\{g\}:g\in G\}=G$, and operations are
given by
\begin{eqnarray}
\{g\}\{h\} & = & \left\{\begin{array}{ll}
\{\bar{g}_-h\}=\{gh_-\} & \mbox{in case $\bar{g}_+=h_+$}  \\
0 & \mbox{in case $\bar{g}_+\neq h_+$}
\end{array}\right. \nonumber \\
\Delta\{g\} & = &  
\sum_{h\bar{k}_+=g,\bar{k}_-=h_-}
\{h\}\otimes\{k\}=
\sum_{h_+k=g,\bar{k}_-=h_-}
\{h\}\otimes\{k\}   \nonumber \\
1 & = & \sum_{g_+\in G_+}\{g_+\} \nonumber \\
\epsilon\{g\} & = &
\left\{\begin{array}{ll}
1 & \mbox{in case $g_+=e$} \\
0 & \mbox{in case $g_+\neq e$} 
\end{array}\right.  \nonumber \\
S\{g\} & = & \{g^{-1}\} \nonumber
\end{eqnarray} 
The basis $G$ is clearly a positive one.
The following theorem shows that, up to rescaling, this is all the finite dimensional
positively based Hopf algebras.

\begin{th}\label{main}
Given any finite dimensional Hopf algebra over ${\bf C}$ with a positive basis
$B$, we can always rescale $B$ by some positive numbers, so that $(H,B)$ is
isomorphic to $(H(G;G_+,G_-),G)$ for a unique group $G$ and a unique factorization
$G=G_+G_-$.
\end{th}

A consequence of the theorem is that the Drinfel'd double $D(H(G;G_+,G_-))$ is also of
the form $H(\tilde{G};\tilde{G}_+,\tilde{G}_-)$ for some group $\tilde{G}$ and
factorization $\tilde{G}=\tilde{G}_+\tilde{G}_-$. By tracing our proof of the theorem,
we find
\[
\begin{array}{ccl}
\tilde{G} & = & G\times G, \\
\tilde{G}_+ & = & \{(g_+,g_-): g_+\in G_+, g_-\in G_-\}\cong G_+\times G_-, \\
\tilde{G}_- & = & \{(g,g): g\in G\}\cong G.
\end{array}
\]

\section{Classification up to Rescaling}

In this section, we prove the classification theorem \ref{main} up to
rescaling. In the next section, we work out the necessary rescaling.

For a positive basis $B$ of $H$, we denote by $B^*=\{b^*:b\in B\}$ the dual basis of
$H^*$. The basis $B^*$ is still positive. We also have the bases $B\otimes B$, 
$B\otimes B^*$, $\cdots$, of $H\otimes H$, 
$H\otimes H^*$, $\cdots$. We say an element $x\in H$ {\em contains} a term $b\in B$ if
the $b$-coordinate of $x$ is not zero. We can say similar things in $H^*$, 
$H\otimes H$, $H\otimes H^*$, $\cdots$.

Given subsets $X$ and $Y$ of a basis $B$, the
notations $X^*$, $X\otimes Y$, $X\otimes Y^*$, etc., carry the obvious meaning as
subsets of $B\otimes B$, $B\otimes B^*$, etc.. The meaning of the phrases like
$X$-term, $X^*$-term, $X\otimes Y$-term, $X\otimes Y^*$-term, etc., are
self-evident.

We will often make use of the duality between $H$ and $H^*$, both being positively
based. This allows us to turn a result about product into a result about
coproduct, and vice versa. The guiding principal for doing this is the following.

\begin{lem}\label{lem0}
The coefficient of $b_1^*\otimes b_2^*$ in $\Delta(b^*)$ is the same as the
coefficient of $b$ in $b_1b_2$. In particular, $\Delta(b^*)$ contains 
$b_1^*\otimes b_2^*$ if and only if $b_1b_2$ contains $b$. Similarly, we may
interchange the role of $H$ and $H^*$.
\end{lem}

\noindent{\em Proof}: Let $b_1b_2=\sum\mu_c c$, and  
$\Delta(b^*)=\sum_{c_1,c_2\in B}\lambda_{c_1,c_2}c_1^*\otimes c_2^*$. Then 
\[
\mu_b=\left<b^*,b_1b_2\right>
=\left<\Delta(b^*),b_1\otimes b_2\right>=\lambda_{b_1,b_2}. 
\]

\hfill$\Box$

After a suitable rescaling by positive numbers, we have
\[
1=d_1+\cdots+d_k, 
\]
where $d_1, \cdots, d_k$ are distinct elements of $B$.

\begin{lem}\label{lem1}
$d_id_j=\delta_{ij}d_i$, and ${\bf C}d_1+\cdots+{\bf C}d_k$ is a
commutative Hopf subalgebra of $H$.
\end{lem}

\noindent{\em Proof}: From $d_i=d_i\cdot 1=d_id_1+\cdots+d_id_k$ and
positivity, $d_id_j$ is a nonnegative multiple of $d_i$.
Similarly, from $d_j=1\cdot d_j$, we see that $d_id_j$ is a nonnegative multiple
of $d_j$. Therefore $d_id_j=0$ when $i\neq j$. Substituting this into
$d_i=d_id_1+\cdots+d_id_k$, we get $d_i=d_id_i$.

From $\Delta (d_1)+\cdots+\Delta (d_k)=\Delta (1)=1\otimes 1=
\sum_{i,j}d_i\otimes d_j$ and positivity, we see that 
${\bf C}d_1+\cdots+{\bf C}d_k$ is closed under $\Delta$.

From $S(d_1)+\cdots+S(d_k)=S(1)=1=d_1+\cdots+d_k$ and positivity, we see that 
${\bf C}d_1+\cdots+{\bf C}d_k$ is closed under $S$.

Therefore ${\bf C}d_1+\cdots+{\bf C}d_k$ is a Hopf subalgebra of $H$.

\hfill$\Box$

Commutative Hopf algebras have been well-classified. We have a group $G_+$, such
that ${\bf C}d_1+\cdots+{\bf C}d_k$ is isomorphic to the dual $({\bf C}G_+)^*$ of
the group algebra. More precisely, we have 
\[
\{d_1, \cdots, d_k\}=G_+^*=\{d_{g_+}: g_+\in G_+\}
\] 
and
\[
d_{g_+}d_{h_+}=\delta_{g_+,h_+}d_{g_+},\hspace{5mm}
\Delta (d_{g_+})=\sum_{h_+k_+=g_+}d_{h_+}\otimes d_{k_+},
\]
\[
S(d_{g_+})=d_{g_+^{-1}},\hspace{5mm}
1=\sum_{g_+\in G_+}d_{g_+},\hspace{5mm}
\epsilon(d_{g_+})=\delta_{g_+,e}
\]

\begin{lem}\label{lem2}
For any $b\in B$, there are unique $g_+, \bar{g}_+\in G_+$, such that
\[
d_{h_+}b=\delta_{g_+,h_+}b,\hspace{5mm}
bd_{h_+}=\delta_{\bar{g}_+,h_+}b.
\]
\end{lem}

\noindent{\em Proof}: From $b=1\cdot b=\sum_{h_+\in G_+}d_{h_+}b$ and
positivity, we have $d_{h_+}b=\lambda_{h_+}b$ for some nonnegative number
$\lambda_{h_+}$. Then from $d_{g_+}d_{h_+}=\delta_{g_+,h_+}d_{g_+}$, we see that 
$\lambda_{g_+}\lambda_{h_+}=\delta_{g_+,h_+}\lambda_{g_+}$. Therefore all except
one $\lambda_*$ vanish. Moreover, the nonvanishing $\lambda_*$ must be 1.

\hfill$\Box$

Denote the unique elements in the last lemma by
\[
\alpha_+(b)=g_+,\hspace{5mm}
\beta_+(b)=\bar{g}_+.
\]
It follows from the lemma that
\begin{equation}\label{a1}
\alpha_+(b)=g_+ 
\Longleftrightarrow d_{g_+}b=b 
\Longleftrightarrow d_{g_+}b\neq 0;
\end{equation}
\begin{equation}\label{a2}
\beta_+(b)=\bar{g}_+ 
\Longleftrightarrow bd_{\bar{g}_+}=b
\Longleftrightarrow bd_{\bar{g}_+}\neq 0.
\end{equation}
We also denote
\begin{eqnarray}
B_{g_+, \bar{g}_+} & = & 
\{b\in B: \alpha_+(b)=g_+,\beta_+(b)=\bar{g}_+\}, \nonumber \\
B_{g_+, \bullet} & = & 
\{b\in B: \alpha_+(b)=g_+\}, \nonumber \\
B_{\bullet, \bar{g}_+} & = &
\{b\in B: \beta_+(b)=\bar{g}_+\}. \nonumber
\end{eqnarray}
Then $B$ is the disjoint union of all $B_{g_+, \bar{g}_+}$.

\begin{lem}\label{lem7}
If $b_1b_2$ contains $b$, then
\begin{equation}\label{x1}
\alpha_+(b_1)=\alpha_+(b),\hspace{5mm}
\beta_+(b_2)=\beta_+(b).
\end{equation}
If $\Delta(b)$ contains $b_1\otimes b_2$, then
\begin{equation}\label{x2}
\alpha_+(b)=\alpha_+(b_1)\alpha_+(b_2),\hspace{5mm}
\beta_+(b)=\beta_+(b_1)\beta_+(b_2).
\end{equation}
\end{lem}

\noindent{\em Proof}: Suppose $b_1b_2$ contains in $b$, and denote
$g=\alpha_+(b)$, $\bar{g}=\beta_+(b)$. Then
$b_1b_2d_{\bar{g}_+}$ and $d_{g_+}b_1b_2$ still contain
$d_{g_+}b=b=bd_{\bar{g}_+}$. Therefore
$b_2d_{\bar{g}_+}\neq 0$ and $d_{g_+}b_1\neq 0$. By (\ref{a1}) and (\ref{a2}), this
implies (\ref{x1}).

Let $\Delta(b)=\sum_{b_1,b_2\in B}\lambda_{b_1,b_2}b_1\otimes b_2$. Then  
\begin{eqnarray}
\Delta(b)=\Delta(d_{g_+}b) & = & \sum_{b_1,b_2\in B,h_+k_+=g_+}
\lambda_{b_1,b_2}d_{h_+}b_1\otimes d_{k_+}b_2 \nonumber \\
& = & \sum_{\alpha_+(b_1)=h_+,\alpha_+(b_2)=k_+,h_+k_+=g_+}
\lambda_{b_1,b_2}b_1\otimes b_2. \nonumber
\end{eqnarray}
Therefore $\Delta(b)$ contains only the terms $b_1\otimes b_2$ satisfying
$\alpha_+(b_1)\alpha_+(b_2)=h_+k_+=g_+=\alpha_+(b)$. The proof for
$\beta_+(b)=\beta_+(b_1)\beta_+(b_2)$ is similar.

\hfill$\Box$

Now let us try to dualize the discussion above.

The dual basis $B^*$ is a
positive basis of the dual Hopf algebra $H^*$. In particular, we have
$\epsilon=\lambda_1e_1^*+\cdots+\lambda_le_l^*$ with $\lambda_i>0$. Clearly, $b^*$
appears in the sum if and only if $\epsilon(b)>0$. Since $d_{e}\in G_+^*$ is
the only element in $G_+^*$ with positive $\epsilon$-value, the intersection
between $\{e_1,\cdots,e_l\}$ and $G_+^*$ is exactly $d_e$. Moreover, by
$1=\epsilon(1)=\epsilon(d_{e})$, we see that the coefficient of $d_{e}^*$ in the
summation for $\epsilon$ is $1$. Therefore we may rescale $e_i$ to get
$\epsilon=e_1^*+\cdots+e_l^*$, without affecting $G_+^*$. 

As in Lemma \ref{lem1}, ${\bf C}e_1^*+\cdots+{\bf C}e_l^*$ is a commutative Hopf
subalgebra of $H^*$. Therefore $\{e_1,\cdots,e_l\}$ form a group
$G_-$ under the multiplication in $H$, with $d_e$ as the unit. To indicate the
group structure, we denote
\[
\{e_1,\cdots,e_l\}=G_-=\{e_{g_-}: g_-\in G_-\}.
\]
Then
\[
e_{g_-}e_{h_-}=e_{g_-h_-},\hspace{5mm}
\Delta (e_{g_-})=e_{g_-}\otimes e_{g_-}+\mbox{terms not in $G_-\otimes G_-$},
\]
\[
S(e_{g_-})=e_{g_-^{-1}}+\mbox{terms not in $G_-$},\hspace{5mm}
\epsilon=\sum_{g_-\in G_-}e_{g_-}^*.
\] 
Since $d_e\in G_+^*$ is also the identity of $G_-$, we will abuse the notation and
also denote $e=d_e$.

The dual of Lemma \ref{lem2} is the following: For any $b^*\in B^*$, we have unique
$g_-,\bar{g}_-\in G_-$, such that
\begin{equation}\label{a'}
e_{h_-}^*b^*=\delta_{\bar{g}_-,h_-}b^*,\hspace{5mm}
b^*e_{h_-}^*=\delta_{g_-,h_-}b^*.
\end{equation}
This enables us to define
\[
\hspace{5mm}
\alpha_-(b)=\bar{g}_-,\hspace{5mm}
\beta_-(b)=g_-.
\]
We may reinterpret (\ref{a'}) as the following lemma. We note that the formula
(\ref{a3}), being equivalent to (\ref{a'}), also serves as the characterizations
of $\alpha_-$, $\beta_-$.

\begin{lem}\label{lem4}
Suppose $\alpha_-(b)=\bar{g}_- $ and $\beta_-(b)=g_-$. Then
\begin{equation}\label{a3}
\Delta(b)=e_{\bar{g}_-}\otimes b+b\otimes e_{g_-}
+\mbox{\em terms not in $B\otimes G_-$ or $G_-\otimes B$}.
\end{equation}
\end{lem}

\noindent{\em Proof}: By Lemma \ref{lem0}, the equalities
$e_{h_-}^*c^*=\delta_{\alpha_-(c),h_-}c^*$ mean that $\Delta(b)$ contains
$e_{h_-}\otimes c$ if and only if $c=b$ and $h_-=\alpha_-(c)$. Moreover, when the
containment happens, the coefficient is 1. This means exactly that 
the only $G_-\otimes B$-term of $\Delta(b)$ is 
$e_{\alpha_-(b)}\otimes b$, and the coefficient is 1. Similarly, the only
$B\otimes G_-$-term of $\Delta(b)$ is $b\otimes e_{\beta_-(b)}$, also with
coefficient 1.

\hfill$\Box$

With the help of Lemma \ref{lem0}, we may also dualize Lemma \ref{lem7}.

\begin{lem}\label{lem7'}
If $b_1\otimes b_2$ is contained in $\Delta(b)$, then
\begin{equation}\label{x1'}
\alpha_-(b_1)=\alpha_-(b),\hspace{5mm}
\beta_-(b_2)=\beta_-(b).
\end{equation}
If $b_1b_2$ contains $b$, then
\begin{equation}\label{x2'}
\alpha_-(b)=\alpha_-(b_1)\alpha_-(b_2),\hspace{5mm}
\beta_-(b)=\beta_-(b_1)\beta_-(b_2).
\end{equation}
\end{lem}

\hfill$\Box$

Having found $G_+$ and $G_-$, we need to further construct $G=G_+G_-$. One approach is
to make use of the Drinfel'd double $D(H)$, which has $B^*\otimes B$ as a positive
basis. From $D(H)$ we can also construct a ``positive group'' $\tilde{G}_+$ and
a ``negative group'' $\tilde{G}_-$. $\tilde{G}_-$ consists of exactly the pairs
$b_1^*\otimes b_2\in B^*\otimes B$, such that 
\[
\epsilon(b_1^*\otimes b_2) = \left< b_1,1 \right> \epsilon(b_2)>0
\]
Thus $\tilde{G}_-=G_+\otimes G_-$ and this is a unique factorization. It remains to
verify that $H(\tilde{G}_-;G_+,G_-)$ is isomorphic to the Hopf algebra $H$.

Instead of making use of Drinfel'd double, we will continue with a more elementary
approach. The approach will also be applicable to the Hopf algebras in the
correspondence category, where the notion of Drinfel'd double is yet to be defined.

\begin{lem}\label{lem5}
The following are equivalent

$(1)$ $\alpha_+(b)=e$; 

$(2)$ $\beta_+(b)=e$;

$(3)$ $b\in G_-$.
 
\noindent In particular, $B_{e,g_+}=B_{g_+,e}=\emptyset$ for $g_+\neq e$.
\end{lem}

\noindent{\em Proof}: Since $d_e=e$ is also the multiplicative unit for $G_-$,
we see (3) implies (1). Conversely, from the last lemma,
$\Delta(b)$ contains $e_{\bar{g}_-}\otimes b$. Then by 
$\epsilon=m(S\otimes 1)\Delta$, $\epsilon(b)$ contains $S(e_{\bar{g}_-})b$, which
further contains $e_{\bar{g}_-^{-1}}b$. If $b$ satisfies (1), then 
$e_{\bar{g}_-}\epsilon(b)$ contains $e_{\bar{g}_-}e_{\bar{g}_-^{-1}}b=d_eb=b$, so
that $\epsilon(b)\neq 0$. This implies $b\sub G_-$. 
Similarly, by using $\epsilon=m(1\otimes S)\Delta$, we can prove (2) implies (3).

\hfill$\Box$

\begin{lem}\label{lem9}
Given $g_+\in G_+$ and $g_-\in G_-$, there is a unique 
$B_{g_+^{-1},\bullet}\otimes B_{g_+,\bullet}$-term $c\otimes b$ contained in
$\Delta(e_{g_-})$. Moreover,  $b$ is the unique element such that
$\alpha_+(b)=g_+$ and $\beta_-(b)=g_-$.
\end{lem}

\noindent{\em Proof}: We begin by finding a 
$B_{g_+^{-1},\bullet}\otimes B_{g_+,\bullet}$-term in $\Delta(e_{g_-})$.
Since $m(1\otimes S)\Delta(e_{g_-})=\epsilon(e_{g_-})=1
=\sum_{h_+\in G_+}d_{h_+}$, $\Delta(e_{g_-})$ must contain a term
$c\otimes b$, such that $cS(b)=\sum_{h_+\in G_+}\lambda_{h_+}d_{h_+}$ with
$\lambda_{g_+^{-1}}>0$. Then we have 
$d_{g_+^{-1}}cS(b)=\lambda_{g_+^{-1}}d_{g_+^{-1}}$. Thus
$d_{g_+^{-1}}c\neq 0$, and by (\ref{a1}), $\alpha_+(c)=g_+^{-1}$. 
Then by Lemmas \ref{lem7}, \ref{lem7'}, and \ref{lem5}, we have
$\alpha_+(b)=\alpha_+(c)^{-1}\alpha_+(e_{g_-})=g_+$.

Next we would like to compare the
$B_{g_+,\bullet}\otimes B_{g_+^{-1},\bullet}\otimes B_{g_+,\bullet}$-terms in
$(\Delta\otimes 1)\Delta(b)$ and $(1\otimes\Delta)\Delta(b)$. For
$(\Delta\otimes 1)\Delta(b)$, such terms can only be obtained from applying 
$\Delta\otimes 1$ to $B\otimes B_{g_+,\bullet}$-terms in $\Delta(b)$. Now for any 
$B\otimes B_{g_+,\bullet}$-term $b_1\otimes b_2$ contained in $\Delta(b)$, we apply
Lemma \ref{lem7} and find
$g_+=\alpha(b)=\alpha_+(b_1)\alpha_+(b_2)=\alpha_+(b_1)g_+$. Therefore
$\alpha_+(b_1)=e$, so that by Lemma \ref{lem5}, we have $b_1\in G_-$. Thus 
$b_1\otimes b_2$ is really a $G_-\otimes B_{g_+,\bullet}$-term contained in
$\Delta(b)$. It then follows from Lemma \ref{lem4} that 
$b_1\otimes b_2=e_{\bar{g}_-}\otimes b$ (where $\bar{g}_-=\alpha_-(b)$), and the 
$B_{g_+,\bullet}\otimes B_{g_+^{-1},\bullet}\otimes B_{g_+,\bullet}$-terms in
$(\Delta\otimes 1)\Delta(b)$ come from $\Delta(e_{\bar{g}_-})\otimes b$.
Similarly, the 
$B_{g_+,\bullet}\otimes B_{g_+^{-1},\bullet}\otimes B_{g_+,\bullet}$-terms in
$(1\otimes \Delta)\Delta(b)$ come from $b\otimes\Delta(e_{g_-})$.

Let $\sum_i\lambda_ic_i\otimes b_i$ be all the 
$B_{g_+^{-1},\bullet}\otimes B_{g_+,\bullet}$-terms in $\Delta(e_{g_-})$. 
Let $\sum_j\lambda'_jb'_j\otimes c'_j$ be all
the $B_{g_+,\bullet}\otimes B_{g_+^{-1},\bullet}$-terms in $\Delta(e_{\bar{g}_-})$. 
Then by comparing the 
$B_{g_+,\bullet}\otimes B_{g_+^{-1},\bullet}\otimes B_{g_+,\bullet}$-terms in 
$\Delta(e_{\bar{g}_-})\otimes b$ and $b\otimes\Delta(e_{g_-})$, we get
\[
\sum_i\lambda_ib\otimes c_i\otimes b_i
=\sum_j\lambda'_jb'_j\otimes c'_j\otimes b.
\]
Consequently, all $b_i=b$, all $b'_j=b$, and
\begin{eqnarray}
\Delta(e_{g_-}) & = & 
(\sum_i\lambda_ic_i)\otimes b 
+\mbox{no $B_{g_+^{-1},\bullet}\otimes B_{g_+,\bullet}$-terms} \label{bbb1} \\
\Delta(e_{\bar{g}_-}) & = & 
b\otimes (\sum_i\lambda_ic_i) 
+\mbox{no $B_{g_+,\bullet}\otimes B_{g_+^{-1},\bullet}$-terms} \label{bbb2}
\end{eqnarray}

The equality (\ref{bbb1}) shows that for any $g_+\in G_+$ and 
$g_-\in G_-$, all the $B_{g_+^{-1},\bullet}\otimes B_{g_+,\bullet}$-terms in 
$\Delta(e_{g_-})$ have the same right component. Now apply this fact to 
$g_+^{-1}\in G_+$ and $\bar{g}_-\in G_-$, we see that all the
$B_{g_+,\bullet}\otimes B_{g_+^{-1},\bullet}$-terms in 
$\Delta(e_{\bar{g}_-})$ have the same right component. Then (\ref{bbb2}) tells us
that $\sum_i\lambda_ic_i$ is really only one term. We thus conclude that for
any $g_+\in G_+$ and $g_-\in G_-$, the 
$B_{g_+^{-1},\bullet}\otimes B_{g_+,\bullet}$-term in $\Delta(e_{g_-})$ is unique.

Now let us prove the uniqueness of the element $b$ satisfying $\alpha_+(b)=g_+$
and $\beta_-(b)=g_-$. We already have one such element appearing in 
\[
\Delta(e_{g_-})
=c\otimes b+\mbox{no $B_{g_+^{-1},\bullet}\otimes B_{g_+,\bullet}$-terms}
\]
Now suppose $b'$ is another base element satisfying $\alpha_+(b')=g_+$
and $\beta_-(b')=g_-$. As before, we write (we do not know apriori that
$\bar{g'}_-=\bar{g}_-$)
\[
\Delta(b')=e_{\bar{g'}_-}\otimes b'+b'\otimes e_{g_-}
+\mbox{terms not in $B\otimes G_-$ or $G_-\otimes B$},
\]
and compare the 
$B_{g_+,\bullet}\otimes B_{g_+^{-1},\bullet}\otimes B_{g_+,\bullet}$-terms in
$(\Delta\otimes 1)\Delta(b')$ and $(1\otimes\Delta)\Delta(b')$. The result is
\[
[\mbox{$B_{g_+,\bullet}\otimes B_{g_+^{-1},\bullet}$-term 
in $\Delta(e_{\bar{g'}_-})$}]\otimes b'
=b'\otimes c\otimes b.
\]
Comparing the rightmost components, we conclude that $b=b'$.

\hfill$\Box$

We make some remarks about the lemma above.

First of all, observe that $c\otimes b$ is also the unique 
$B\otimes B_{g_+,\bullet}$-term (as well as the unique
$B_{g_+^{-1},\bullet}\otimes B$-term) in $\Delta(e_{g_-})$. As a matter of fact, by
Lemma \ref{lem7}, any $B\otimes B_{g_+,\bullet}$-term in $\Delta(e_{g_-})$ must be a
$B_{g_+^{-1},\bullet}\otimes B_{g_+,\bullet}$-term.

Secondly, if we start with $\bar{g}_+\in G_+$ and $\bar{g}_-\in G_-$, similar
proof also tells us that there is a unique 
$B_{\bullet,\bar{g}_+}\otimes B_{\bullet,\bar{g}_+^{-1}}$-term $b\otimes c$ contained
in $\Delta(e_{\bar{g}_-})$. Moreover, $b$ is the unique element such that
$\beta_+(b)=\bar{g}_+$ and $\alpha_-(b)=\bar{g}_-$. Briefly speaking, one may
start with $m(S\otimes 1)\Delta(e_{\bar{g}_-})=1$ and conclude
$\Delta(e_{\bar{g}_-})$ contains a $B_{\bullet,\bar{g}_+}\otimes
B_{\bullet,\bar{g}_+^{-1}}$-term. Then by comparing the 
$B_{\bullet,\bar{g}_+}\otimes B_{\bullet,\bar{g}_+^{-1}}\otimes
B_{\bullet,\bar{g}_+}$-terms in 
$(\Delta\otimes 1)\Delta(b)$ and $(1\otimes\Delta)\Delta(b)$, our claim
follows. In particular, we also see that $b\ra (\beta_+(b),\alpha_-(b))$ is a
one-to-one correspondence.

Thirdly, we have the one-to-one correspondences
\begin{equation}\label{eee}
(g_+,g_-)
\stackrel{(\alpha_+,\beta_-)}{\longleftarrow}
b
\stackrel{(\beta_+,\alpha_-)}{\longrightarrow}
(\bar{g}_+,\bar{g}_-),
\end{equation}
and the equalities
\begin{eqnarray}
\Delta(e_{g_-}) & = & \lambda c\otimes b+\cdots \label{eee1} \\
\Delta(e_{\bar{g}_-}) & = & \lambda b\otimes c+\cdots \label{eee2}
\end{eqnarray}
The one-to-one correspondences (\ref{eee}) enables us to define
a right action of $G_-$ on $G_+$ 
\begin{equation}\label{action1}
G_+\times G_-\ra G_+: \hspace{5mm}
(g_+,g_-)\mapsto g_+^{\;g_-}=\bar{g}_+,
\end{equation}
and a left action of $G_+$ on $G_-$ 
\begin{equation}\label{action2}
G_+\times G_-\ra G_-: \hspace{5mm}
(g_+,g_-)\mapsto \;^{g_+}g_-=\bar{g}_-.
\end{equation}
These are presumably the actions (\ref{gpd8}) and (\ref{gpd9}) in our standard model.
In Lemma \ref{lem13}, we will show they are indeed actions and satisfy matching
relations (\ref{gpdmatch}).

The action (\ref{action1}) has the following interpretation: Given $g_+$ and
$g_-$, we use (\ref{eee1}) to find the unique $b\in B_{g_+,\bullet}$. Then
$\bar{g}_+=\beta_+(b)$ is $\;^{g_-}g_+$ (i.e., $b\in B_{g_+,\bar{g}_+}$).

The action (\ref{action2}) has the following interpretation: Given $g_+$ and
$g_-$, we establish (\ref{eee1}) and (\ref{eee2}) with $b\in B_{g_+,\bullet}$. Then
the element $\bar{g}_-$ in (\ref{eee2}) is $g_-^{\;g_+}$.

The reason for the interpretations above follows from (\ref{bbb1}),
(\ref{bbb2}), and the uniqueness. We may also find $g_+$, $g_-$ from
$(\bar{g}_+,\bar{g}_-)$ in the similar way. These are the other two actions in the
standard model, and we will not use them.

\begin{lem}\label{lem11}
For any $b\in B_{g_+,\bar{g}_+}$, there is $b'\in B_{\bar{g}_+,g_+}$, such that 
$bb'=\lambda d_{g_+}$ and $b'b=\lambda d_{\bar{g}_+}$ for some $\lambda>0$.
\end{lem}

\noindent{\em Proof}: Let $\beta_-(b)=g_-$. From Lemma \ref{lem9} and the first
two of the subsequent remarks, we have a unique $B\otimes B_{\bullet,\bar{g}_+}$-term
$c\otimes b$ in $\Delta(e_{g_-})$. Then by 
$m(S\otimes 1)\Delta(e_{g_-})=\epsilon(e_{g_-})=1$, we have $S(c)b=\mu
d_{\bar{g}_+}$ for some $\mu>0$. By positivity, there is some term $b'$ in $S(c)$,
such that
$b'b=\lambda d_{\bar{g}_+}$ for some
$\lambda>0$. 

We claim that $b'\in B_{\bar{g}_+,g_+}$. From $d_{\bar{g}_+}b'b=\lambda
d_{\bar{g}_+}^2=\lambda d_{\bar{g}_+}\neq 0$, we have $d_{\bar{g}_+}b'\neq 0$, so
that $\alpha_+(b')=g_+$. From $(b'd_{g_+})b=b'(d_{g_+}b)=b'b=\lambda
d_{\bar{g}_+}\neq 0$, we have $b'd_{g_+}\neq 0$, so that
$\beta_+(b')=\bar{g}_+$.

We have shown that for any $b\in B_{g_+,\bar{g}_+}$, there is
$b'\in B_{\bar{g}_+,g_+}$, such that $b'b=\lambda d_{\bar{g}_+}$ for some
$\lambda>0$. Applying this conclusion to $b'$, there is $b''\in
B_{g_+,\bar{g}_+}$, such that $b''b'=\lambda' d_{g_+}$ for some
$\lambda'>0$. Then we have
\[
b''b'b=(b''b')b=\lambda' d_{g_+}b=\lambda' b,\hspace{5mm}
b''b'b=b''(b'b)=\lambda b''d_{\bar{g}_+}=\lambda b''.
\]
This implies that $\lambda=\lambda'$ and $b''=b$. 

\hfill$\Box$

The next lemma shows that the antipode $S$ is a permutation on $B$ up to rescaling.
We give a proof that is also valid in the correspondence category. For Hopf
algebras over ${\bf C}$, it can also be proved more directly by using the
positivity of $S$ and the fact that
$S^N=id$ for some positive integer $N$.

\begin{lem}\label{lem6}
Given any $b\in B$, $S(b)=\lambda b'$ for some $b'\in B$ and $\lambda>0$. Moreover,
we have $\alpha_-(b')=\beta_-(b)^{-1}$, $\beta_+(b')=\alpha_+(b)^{-1}$, and these
uniquely determine $b'$.
\end{lem}

\noindent{\em Proof}: First we show that $S(e_{g_-})=e_{g_-^{-1}}$. We already
know $S(e_{g_-})=e_{g_-^{-1}}+\mbox{terms not in $G_-$}$. Suppose $g_-\neq e$, and
$S(e_{g_-})$ contains $b\not\in G_-$. Then by Lemma \ref{lem5},
$h_+=\alpha_+(b)\neq e$, so that
$d_{h_+}S(e_{g_-})=S(e_{g_-}d_{h_+^{-1}})=0$. On the other hand,
$d_{h_+}S(e_{g_-})$ still contains $d_{h_+}b=b$ and cannot vanish. The
contradiction shows that $S(e_{g_-})=e_{g_-^{-1}}$. On the other hand, suppose
$g_-=e$. Then we already know $e_{g_-}=d_e$ and $S(d_e)=d_e$.

Let $g_{\pm}$ and $\bar{g}_{\pm}$ be associated to
$b$ as usual. Then from Lemma \ref{lem9} we have
\[
\Delta(e_{g_-})=\lambda c\otimes b+\cdots
\]
Applying $S$, we have
\[
\Delta(e_{g_-^{-1}})=\lambda S(b)\otimes S(c)+\cdots
\]
The proof of the last lemma implies that $S(c)$ contains a term 
$c'\in B_{g_+,\bar{g}_+}$. Therefore if $S(b)$ contains $b'$, then $b'\otimes c'$
must be contained in $\Delta(e_{g_-^{-1}})$. By the first remark after Lemma
\ref{lem9}, the $B\otimes B_{g_+,\bar{g}_+}$-term in $\Delta(e_{g_-^{-1}})$ must
be unique. Therefore for the fixed $c'$, there is only one term like 
$b'\otimes c'$ in $\Delta(e_{g_-^{-1}})$. This implies that there is only one term, a
constant multiple of $b'$, contained in $S(b)$. 

Thus we conclude the first part of the lemma, and have 
\[
S(b)=\mu b',\hspace{3mm}
S(c)=\nu c',\hspace{5mm}
\Delta(e_{g_-^{-1}})=\lambda\mu\nu b'\otimes c'+\cdots.
\]
Moreover, we have also assumed $c'\in B_{g_+,\bar{g}_+}$.
By lemmas \ref{lem7} and \ref{lem7'}, we then conclude that 
$\alpha_-(b')=\beta_-(b)^{-1}$ and
$\beta_+(b')=\alpha_+(b)^{-1}$. By the second remark after Lemma \ref{lem9}, this
uniquely determines $b'$.

\hfill$\Box$

\begin{lem}\label{lem12}
For $b\in B_{g_+,\bar{g}_+}$ and $c\in B_{h_+,\bar{h}_+}$, there are two
possibilities for $bc$:

$(1)$ If $\bar{g}_+=h_+$, then $bc$ is a positive multiple of an element in
$B_{g_+,\bar{h}_+}$;

$(2)$ If $\bar{g}_+\neq h_+$, then $bc=0$.
\end{lem}

\noindent{\em Proof}: If $\bar{g}_+\neq h_+$, then
$bc=(bd_{\bar{g}_+})(d_{h_+}c)=b(d_{\bar{g}_+}d_{h_+})c=0$.

If $\bar{g}_+=h_+$, then Lemma \ref{lem11} immediately
implies $bc\neq 0$. Suppose $bc$ contains two different terms $d, d'\in B$.
Then by $d_{g_+}bc=bc=bcd_{\bar{h}_+}$ and positivity, we see that 
$d,d'\in B_{g_+,\bar{g'}_+}$. By Lemma \ref{lem11}, there is $b'\in
B_{\bar{g}_+,g_+}$, such that $bb'=\lambda d_{g_+}$ and $b'b=\lambda d_{\bar{g}_+}$
for some $\lambda>0$. Then $b'bc=\lambda c$ contains $b'd$ and $b'd'$. By
positivity, we must have $b'd=\mu c$ and 
$b'd'=\mu' c$ for some $\mu>0$ and $\mu'>0$. Thus we
conclude that
\[
\lambda d=\lambda d_{g_+}d=bb'd=\mu bc,\hspace{5mm}
\lambda d'=\lambda d_{g_+}d'=bb'd'=\mu' bc.
\]
This is in contradiction with the assumption that $d$ and $d'$ are distinct.

\hfill$\Box$

Lemma \ref{lem6} implies that $S^{-1}$ exists, and for any $b\in B$, the coordinates of
$S^{-1}(b)$ are non-negative. Therefore $(H^*)^{op}$ is also positively
based. In particular, all the proofs we have done so far are valid in $(H^*)^{op}$.
Note that the two actions (\ref{action1}) and (\ref{action2}) are
related by duality in the following way.
\[
\begin{array}{ccccc}
H & \Longleftrightarrow & H^* & \Longleftrightarrow & (H^*)^{op} \\
\alpha_+(b)=g_+ & \longleftrightarrow & 
\alpha_-(b^*)=\bar{g}_- & \longleftrightarrow &
\alpha_-^{op}(b^*)=g_-   \\
\beta_+(b)=\bar{g}_+ & \longleftrightarrow & 
\beta_-(b^*)=g_- & \longleftrightarrow &
\beta_-^{op}(b^*)=\bar{g}_-   \\
\alpha_-(b)=\bar{g}_- & \longleftrightarrow & 
\alpha_+(b^*)=g_+ & \longleftrightarrow &
\alpha_+^{op}(b^*)=\bar{g}_+   \\
\beta_-(b)=g_- & \longleftrightarrow & 
\beta_+(b^*)=\bar{g}_+ & \longleftrightarrow &
\beta_+^{op}(b^*)=g_+   \\
\bar{g}_+=g_+^{\;g_-} & \longleftrightarrow &
g_-=\bar{g}_-^{\;\bar{g}_+} & \longleftrightarrow &
\bar{g}_-=\;^{g_+}g_- \\
\bar{g}_-=\;^{g_+}g_- & \longleftrightarrow &
g_+=\;^{\bar{g}_-}\bar{g}_+ & \longleftrightarrow &
\bar{g}_+=g_+^{\;g_-}
\end{array}
\] 
Therefore by restating the results we have proved in $(H^*)^{op}$,
we may exchange the two actions in these results.

\begin{lem}\label{lem13}
$(\ref{action1})$ and $(\ref{action2})$ are indeed actions. Moreover, the two actions
satisfy the matching relations $(\ref{gpdmatch})$.
\end{lem}

\noindent{\em Proof}: Given $g_+^{\;g_-}=h_+$ and $h_+^{\;h_-}=k_+$, we may find
(by the third remark after Lemma \ref{lem9}) 
$b\in B_{g_+,h_+}$ and $b'\in B_{h_+,k_+}$, with $g_-=\beta_-(b)$ and
$h_-=\beta_-(b')$. Then by Lemma \ref{lem12}, we have 
$bb'=\lambda c$, with $\lambda>0$ and $c\in B_{g_+,k_+}$. 
By Lemma \ref{lem7'}, we have $\beta_-(c)=g_-h_-$. Therefore 
$g_+^{\;g_-h_-}=k_+=h_+^{\;h_-}=(g_+^{\;g_-})^{\;h_-}$. This proves
(\ref{action1}) is a right action. 

As explained in the remark above, we may exchange the two
actions and still get an equality. Applying the principle to
$g_+^{\;g_-h_-}=(g_+^{\;g_-})^{h_-}$, we get
$\;^{\;h_+g_+}g_-=\;^{\;h_+}(\;^{g_+}g_-)$. This proves (\ref{action2}) is a
left action. 

Now we prove the matching relations.
By the third remark after Lemma \ref{lem9}, we interpret the actions
$\bar{g}_+=g_+^{\;g_-}$ and $\bar{g}_-=\;^{g_+}g_-$ as the following equalities
\begin{equation}\label{c1}
\left\{
\begin{array}{ccl}
\Delta(e_{g_-}) & = & \lambda c\otimes b+\cdots,  \\
\Delta(e_{\bar{g}_-}) & = & \lambda b\otimes c+\cdots, 
\end{array}
\right.
\hspace{5mm}
b\in B_{g_+,\bar{g}_+}.
\end{equation}
We also interpret the action
$\bar{h}_-=\;^{\bar{g}_+}h_-$ as the following equalities
\begin{equation}\label{c2}
\left\{
\begin{array}{ccl}
\Delta(e_{h_-}) & = & \lambda' c'\otimes b'+\cdots,  \\
\Delta(e_{\bar{h}_-}) & = & \lambda' b'\otimes c'+\cdots, 
\end{array}
\right.
\hspace{5mm}
b\in B_{\bar{g}_+,\bullet}.
\end{equation}
Multiplying (\ref{c1}) and (\ref{c2}) together, we have (where we use Lemma
\ref{lem12} so that $cc'$ and $bb'$ are elements of $B$ up positive rescaling)
\begin{equation}\label{c3}
\left\{
\begin{array}{ccl}
\Delta(e_{g_-h_-})=\Delta(e_{g_-})\Delta(e_{h_-}) 
& = & \lambda\lambda' cc'\otimes bb'+\cdots,  \\
\Delta(e_{\bar{g}_-\bar{h}_-})=\Delta(e_{\bar{g}_-})\Delta(e_{\bar{h}_-}) 
& = & \lambda\lambda' bb'\otimes cc'+\cdots, 
\end{array}
\right.
\hspace{5mm}
bb'\in B_{g_+,\bullet}.
\end{equation}
By the third remark after Lemma \ref{lem9}, (\ref{c3}) means
$\bar{g}_-\bar{h}_-=\;^{g_+}(g_-h_-)$. Therefore
\[
\;^{g_+}(g_-h_-)=\bar{g}_-\bar{h}_-=\;^{g_+}g_-\;^{\bar{g}_+}h_-
=\;^{g_+}g_-\;\;^{(g_+^{\;g_-})}h_-.
\]
This is the matching relation. Applying remark made before the lemma, we get the
other matching relation
\[
(h_+g_+)^{g_-}=h_+^{(\;^{g_+}g_-)}\;g_+^{\;g_-}.
\]

\hfill$\Box$

With a matching group pair $G_+, G_-$, we may construct a group $G=G_+G_-$ and
projections $\alpha_{\pm},\beta_{\pm}:G\ra G_{\pm}$. Furthermore, we have a
standard Hopf algebra $H(G;G_+,G_-)$ with $G$ as a positive basis. The next Lemma
shows that up to positive scaling, the product between the basis elements of $H$
is the same as the product between the basis elements of $H(G;G_+,G_-)$.

\begin{lem}\label{lem10}
$bb'$ is a positive multiple of $c$ if and only if
\[
\alpha_+(b)=\alpha_+(c),\hspace{5mm}
\beta_+(b)=\alpha_+(b'),\hspace{5mm}
\alpha_-(b)\alpha_-(b')=\alpha_-(c).
\]
Moreover, for fixed $c$, $(b,b')\mapsto \alpha_-(b)$ is a one-to-one
correspondence between the pairs $(b,b')$ and $G_-$.
\end{lem}

\noindent{\em Proof}: The necessity follows from (\ref{x1}), (\ref{x2'}), and
Lemma \ref{lem12}. For sufficiency, the condition $\beta_+(b)=\alpha_+(b')$ implies
$bb'=\lambda c'$ for some $\lambda>0$ and $c'\in B$. Then we have
$\alpha_+(c)=\alpha_+(b)=\alpha_+(c')$ and
$\alpha_-(c)=\alpha_-(b)\alpha_-(b')=\alpha_-(c')$. If we can show
$(\alpha_+(c),\alpha_-(c))$ are in one-to-one correspondence with
$(\alpha_+(c),\beta_-(c))$, then it follows from Lemma
\ref{lem9} that $c=c'$.

By the definition of the action of $G_+$ on $G_-$, from
$(\alpha_+(c),\beta_-(c))$ we have $\alpha_-(c)=\;^{\alpha_+(c)}\beta_-(c)$.
Since this is a group action, from $(\alpha_+(c),\beta_-(c))$ we also have
$\beta_-(c)=\;^{\alpha_+(c)^{-1}}\alpha_-(c)$.
This proves the sufficiency.

For the one-to-one correspondence, we fix $\alpha_-(b)\in G_-$. The condition
$\alpha_+(b)=\alpha_+(c)$ then uniquely determines $b$. Moreover, the conditions
$\beta_+(b)=\alpha_+(b')$, $\alpha_-(b)\alpha_-(b')=\alpha_-(c)$ give us
$\alpha_+(b')$ and $\alpha_-(b')$, which also has a unique corresponding $b'$.
This proves the one-to-one correspondence between pairs $(b,b')$ and $G_-$.

\hfill$\Box$

By making use of Lemma \ref{lem0} and the remark before Lemma \ref{lem13},
we get the following dual of Lemma \ref{lem10}. It shows that up
to positive scaling, the coproduct of basis elements of $H$ is the same as the
coproduct of basis elements of $H(G;G_+,G_-)$.

\begin{lem}\label{lem10'}
For any $b\in B$, $\Delta(b)$ contains exactly those $c\otimes c'$ satisfying
\[
\alpha_-(c)=\alpha_-(b),\hspace{5mm}
\beta_-(c)=\alpha_-(c'),\hspace{5mm}
\alpha_+(c)\alpha_+(c')=\alpha_+(b).
\]
Moreover, for fixed $b$, $c\otimes c'\mapsto \alpha_+(c)$ is a one-to-one
correspondence between $c\otimes c'$ and $G_+$.
\end{lem}

\hfill$\Box$

\section{Rescaling}

Given a finite dimensional Hopf algebra $H$ with a positive basis $B$, we have
identified $H$ with a standard Hopf algebra $H(G;G_+,G_-)$ (as well as $B$ with
$G$) up to positive scaling. In this section, we try to positively rescale the
elements of $G$, so that the Hopf algebra structure in
$H$ and $H(G;G_+,G_-)$ are exactly the same.

The problem may be viewed in the following way. On the vector space $V={\bf C}G$,
we have the standard Hopf algebra structure $H(G;G_+,G_-)$. We also have the given
Hopf algebra structure $H$. The
problem is to find a function $\tau:G\ra {\bf R}_+$, such that 
$g\mapsto \tau(g)g$ is a Hopf algebra isomorphism from $H$ to $H(G;G_+,G_-)$.

The Hopf algebra structure $H$ is determined by 
\begin{eqnarray}
\{g\}\{h\} & = & 
\lambda(g,h)\{gh_-\} \hspace{3mm}\mbox{in case $\bar{g}_+=h_+$} \nonumber \\
\Delta\{g\} & = & 
\sum_{h_+k=g,\bar{k}_-=h_-}\mu(h,k)\{h\}\otimes\{k\} \label{s1} \\
S\{g\} & = & \nu(g)\{g^{-1}\}  \nonumber
\end{eqnarray}
The structure constants $\lambda(g,h)$, $\mu(h,k)$, $\nu(g)$ are positive and
defined exactly when the conditions indicated above are satisfied. We assume $H$
is already well-scaled on $G_+^*$ and $G_-$, i.e., with canonical products,
coproducts, and antipodes, $({\bf C}G_+)^*$ is a sub Hopf
algebra of $H$ and $({\bf C}G_-)^*$ is a sub Hopf algebras of $H^*$. Then the
Hopf algebra condition means the following compatibility equations (in cases both
sides are defined):
\begin{eqnarray}
\lambda(g,h)\lambda(gh_-,k) & = & \lambda(g,hk_-)\lambda(h,k), \nonumber \\
\mu(g,h)\mu(g_+h,k) & = & \mu(g,h_+k)\mu(h,k),  \nonumber \\
\lambda(h_+k,h'_+k')\mu(hh'_-,kk'_-) & = & 
\lambda(h,h')\lambda(k,k')\mu(h,k)\mu(h',k'), \nonumber  \\
\lambda(h_+,h) & = & 1, \nonumber  \\
\lambda(h,\bar{h}_+) & = & 1, \nonumber \\
\mu(h,h_-) & = & 1, \nonumber \\
\mu(\bar{h}_-,h) & = & 1, \nonumber \\
\nu(h)\lambda(h^{-1},h_+^{-1}\bar{h}_-) & = & 1, \nonumber \\
\nu(h)\lambda(h_+^{-1}\bar{h}_-,h) & = & 1. \nonumber
\end{eqnarray}
It is clear from the equations above that for any
real number $\epsilon$, $\lambda^{\epsilon}$, $\mu^{\epsilon}$, $\nu^{\epsilon}$
still satisfy the compatibility equations. Therefore if we replace $\lambda$,
$\mu$, $\nu$ in (\ref{s1}) with $\lambda^{\epsilon}$, $\mu^{\epsilon}$,
$\nu^{\epsilon}$, we still get a Hopf algebra (which is also well-scaled on
$G_+^*$ and $G_-$), which we denote by $H_{\epsilon}$. Then we have $H=H_1$ and
$H(G;G_+,G_-)=H_0$.

\begin{lem}
The Hopf algebras $H_{\epsilon}$ are semisimple and cosemisimple.
\end{lem}

\noindent{\em Proof}: Let $\Lambda=\sum_{g_-\in G_-}\{g_-\}$. Because
$H_{\epsilon}$ is already well-scaled on $G_-$, we have
$x\Lambda=\epsilon(x)\Lambda$. Therefore $\Lambda$ is a left integral of the Hopf
algebra. Since $\epsilon(\Lambda)=|G_-|\neq 0$, it follows from \cite{LS} that
$H_{\epsilon}$ is semi-simple.

Similarly, the integral $\sum_{g_+\in G_+}\{g_+\}^*$ tells us that
$H_{\epsilon}^*$ is also semi-simple.

\hfill$\Box$

Consider $H_0$ as a point in the space $X$ of all bialgebra structures on $V$. We
have natural action of the general linear group $GL(V)$ on $X$. The orbit
$GL(V)H_0$ consists of all the bialgebra structures linearly isomorphic to $H_0$.
According to the Corollary 1.5 and the Theorem 2.1 of \cite{S}, the semisimplicity
and the cosemisimplicity of $H_{\epsilon}$ implies that the orbit $GL(V)H_0$ is
a Zariski open subset of $X$. Therefore for sufficiently small
$\epsilon$, $H_{\epsilon}$ belongs to the orbit $GL(V)H_0$. In particular, for any
$\delta>0$, we can find an $\epsilon>0$ and a linear transformation $T$, such
that $T(H_{\epsilon})=H_0$ and 
\begin{equation}\label{s10}
T\{g\}=\sum_{g\in G}\tau_{g,h}\{h\}, \hspace{3mm}
|\tau_{g,g}-1|<\delta,\hspace{3mm}
|\tau_{g,h}|<\delta \hspace{3mm} \mbox{for $h\neq g$}.
\end{equation}

\begin{lem}\label{lem50}
For sufficiently small $\epsilon$, $T$ preserves ${\bf C}G_-$.
\end{lem}

\noindent{\em Proof}: 
Consider $T\{g_-\}=\sum_{h\in G}\tau_{g_-,h}\{h\}$. We would like to show
$\tau_{g_-,h}=0$ for $h\not\in G_-$. If not, we find a term $\tau_{g_-,k}\{k\}$,
$k\not\in G_-$, with
\[
|\tau_{g_-,k}|=\mbox{max}\{|\tau_{g_-,h}|: h\not\in G_-\}.
\]
Since $k\not\in G_-$, we have $k_+\neq e$. Therefore $\{k_+\}\{g_-\}=0$ in
$H_{\epsilon}$. Since $T$ preserves products, we have $T\{k_+\}T\{g_-\}=0$.
On the other hand, we will show that for sufficiently small $\epsilon$, the
coefficient of $\{k\}$ in
$T\{k_+\}T\{g_-\}$ cannot be zero. The contradiction proves the Lemma.

Write
\begin{eqnarray}
T\{k_+\}
& = & \tau_{k_+,k_+}\{k_+\}+\sum_{h\neq k_+}\tau_{k_+,h}\{h\}
=P+Q; \nonumber \\ 
T\{g_-\}
& = & \tau_{g_-,k}\{k\}+\sum_{h\in G_-}\tau_{g_-,h}\{h\}
+\sum_{h\not\in G_-,h\neq k}\tau_{g_-,h}\{h\}
= L+M+N. \nonumber
\end{eqnarray}
Suppose (\ref{s10}) is satisfied for some small $\delta>0$. Then we have
\begin{enumerate}
\item $PL=\tau_{k_+,k_+}\tau_{g_-,k}\{k\}$, with
$|\tau_{k_+,k_+}\tau_{g_-,k}|\geq(1-\delta)|\tau_{g_-,k}|$;
\item $P(M+N)$ can never produce $\{k\}$;
\item $Q(L+N)$ may produce some $\{k\}$-terms. However, the
total number of product terms is no more than $(\dim H)^2$, and the coefficients
are all smaller than $\delta|\tau_{g_-,k}|$. Therefore the $k$-coefficient of this
part is bounded by $\delta|\tau_{g_-,k}|(\dim H)^2$;
\item $QM$ is a sum of $G_-$-terms, so that it can never produce $\{k\}$.
\end{enumerate}
Now for sufficiently small $\epsilon$, we can find $T$ satisfying (\ref{s10}),
with $\delta$ small enough so that $(1-\delta)>\delta(\dim H)^2$. Then we conclude
that the $\{k\}$-term in $T\{k_+\}T\{g_-\}$ is nontrivial. 

\hfill$\Box$

\begin{lem}\label{lem51}
For sufficiently small $\epsilon$, $T({\bf C}(G_+-\{e\}))\sub {\bf C}(G-G_-)$.
\end{lem}

\noindent{\em Proof}: For any $g_+\in G_+$, $g_+\neq e$, write
\[
T\{g_+\}
=\sum_{h_-\in G_-}\tau_{g_+,h_-}\{h_-\}
+\sum_{h\not\in G_-}\tau_{g_+,h}\{h\}
=P+Q.
\]
If $P$ is nontrivial, then we have some $k_-\in G_-$, such that
\[
|\tau_{g_+,k_-}|=\mbox{max}\{|\tau_{g_+,h_-}|: h_-\in G_-\}.
\]
Since $g_+\neq e$, we have $\{g_+\}\{e\}=0$ in $H_{\epsilon}$, which implies
$T\{g_+\}T\{e\}=0$. On the other hand, we will show that for sufficiently small
$\epsilon$, the coefficient of $\{k_-\}$ in $T\{g_+\}T\{e\}$ cannot be zero.
The contradiction proves the lemma.

For sufficiently small $\epsilon$, we have 
\[
T\{e\}=\tau_{e,e}\{e\}+\sum_{h_-\in G_-,h_-\neq e}\tau_{e,h_-}\{h_-\}=L+M
\]
from Lemma \ref{lem50}. Suppose (\ref{s10}) is further satisfied for some small
$\delta>0$. Then we have
\begin{enumerate}
\item The $k_-$-term in $PL$ is $\tau_{g_+,k_-}\tau_{e,e}\{k_-\}$, with
$|\tau_{g_+,k_-}\tau_{e,e}|>(1-\delta)|\tau_{g_+,k_-}|$;
\item $PM$ may produce some $\{k\}$-terms. However, the
total number of product terms is no more than $|G_-|^2$, and the coefficients are
all smaller than $\delta|\tau_{g_+,k_-}|$. Therefore $k$-coefficient of this
part is bounded by $\delta|\tau_{g_+,k_-}||G_-|^2$;
\item $Q(L+M)$ is contains no $G_-$-terms, so that it can never produce $\{k_-\}$.
\end{enumerate}
Now for sufficiently small $\epsilon$, we can find $T$ satisfying Lemma \ref{lem50}
and (\ref{s10}), with $\delta$ small enough so that $(1-\delta)>\delta|G_-|^2$.
Then we conclude that the $\{k_-\}$-term in $T\{g_+\}T\{e\}$ is nontrivial.

\hfill$\Box$

\begin{lem}\label{lem52}
For sufficiently small $\epsilon$, $T$ preserves ${\bf C}(G-G_-)$.
\end{lem}

\noindent{\em Proof}: Let $g\not\in G_-$. Then $g_+\neq e$, so that by 
Lemma \ref{lem51}, $T\{g_+\}$ contains no $G_-$-terms.

In $H_{\epsilon}$, we have $\{g_+\}\{g\}=\{g\}$, which implies
$T\{g\}=T\{g_+\}T\{g\}$. Since $T\{g_+\}$ contains no $G_-$-terms,
$T\{g_+\}T\{g\}$ also contains no $G_-$-terms. Consequently, $T\{g\}$ contains no
$G_-$-terms.

\hfill$\Box$

The Lemmas \ref{lem50} and \ref{lem52} shows that if we decompose ${\bf C}G$ into 
${\bf C}G_-\oplus {\bf C}(G-G_-)$, then $T:H_{\epsilon}\ra H_0$ takes the form
\begin{equation}\label{s11}
T=\left(\begin{array}{cc} * & 0 \\ 0 & * \end{array}\right)
\begin{array}{l} \leftarrow G_- \\ \leftarrow G-G_-\end{array}.
\end{equation}
Now we apply this to the dual $(T^*)^{-1}:H_{\epsilon}^*\ra H_0^*$. Of course to
do so may requires us to choose even smaller $\epsilon$, so that all the lemmas
work for $(T^*)^{-1}$. Since the groups $G_+$ and $G_-$ derived from $H^*$ are the
same as the groups $G_-$ and $G_+$ derived from $H$, we conclude that with regard
to the decomposition 
$({\bf C}G)^*=({\bf C}G_+)^*\oplus ({\bf C}(G-G_+))^*$,
\[
(T^*)^{-1}=\left(\begin{array}{cc} * & 0 \\ 0 & * \end{array}\right)
\begin{array}{l} \leftarrow G_+^* \\ \leftarrow G^*-G_+^*\end{array}.
\]
Consequently, we get
\[
T=\left(\begin{array}{cc} * & 0 \\ 0 & * \end{array}\right)
\begin{array}{l} \leftarrow G_+ \\ \leftarrow G-G_+\end{array}.
\]
In particular, $T$ preserves ${\bf C}G_+$.

\begin{lem}\label{lem54}
For sufficiently small $\epsilon$, $T$ is identity on $G_+$.
\end{lem}

\noindent{\em Proof}: Let $g_+\in G_+$. For sufficiently small $\epsilon$,
from the discussion above we have
\[
T\{g_+\}=\sum_{h_+\in G_+}\tau_{g_+,h_+}\{h_+\},
\]
where $\tau_{g_+,g_+}$ is close to 1 and $\tau_{g_+,h_+}$ is small for $h_+\neq
g_+$. Moreover, since $H_{\epsilon}$ is already well-scaled on $G_+^*$, we have
$\{g_+\}^2=\{g_+\}$ in $H_{\epsilon}$. Therefore in $H_0$, we have
\[
\sum_{h_+\in G_+}\tau_{g_+,h_+}\{h_+\}
=T\{g_+\}=(T\{g_+\})^2=\sum_{h_+\in G_+}\tau_{g_+,h_+}^2\{h_+\}.
\]
Comparing the two sides, we have $\tau_{g_+,h_+}=\tau_{g_+,h_+}^2$. Since
$\tau_{g_+,g_+}$ is the only coefficient close to 1, we
conclude that $\tau_{g_+,h_+}=\delta_{g_+,h_+}$.

\hfill$\Box$

The lemma also applies to the dual Hopf algebras. Taking the adjoint and then the
inverse of (\ref{s11}), we have
\begin{equation}\label{s14}
(T^*)^{-1}=\left(\begin{array}{cc} * & 0 \\ 0 & * \end{array}\right)
\begin{array}{l} \leftarrow G_-^* \\ \leftarrow G^*-G_-^*\end{array}.
\end{equation}
Since $G_-^*$ in $H^*$ plays the role of $G_+$ in $H$, the dual of Lemma
\ref{lem54} then shows that $(T^*)^{-1}$ is identity on $G_-^*$, so that
\begin{equation}\label{s15}
(T^*)^{-1}=\left(\begin{array}{cc} I & 0 \\ 0 & * \end{array}\right)
\begin{array}{l} \leftarrow G_-^* \\ \leftarrow G^*-G_-^*\end{array}.
\end{equation}
By taking the inverse and then the adjoint of (\ref{s15}), we see that $T$ is also
identity on $G_-$.

Having proved $T$ is identity on $G_{\pm}$, we are ready to show $T$ is a
scaling everywhere.

\begin{lem}\label{lem55}
For sufficiently small $\epsilon$, $T$ is a scaling on $G$.
\end{lem}

\noindent{\em Proof}: For any $g=g_+g_-\in G$, we have $\{g_+\}\{g\}=\{g\}$ in
$H_{\epsilon}$ and $T\{g_+\}=\{g_+\}$. Then
$T\{g\}=T(\{g_+\}\{g\})=\{g_+\}T\{g\}$ implies that
$T\{g\}$ is a linear combination of terms $\{h\}$ with $h_+=g_+$. This is
equivalent to the fact that, with regard to the disjoint union decomposition
$\coprod_{g_+\in G_+}\{h\in G:\alpha_+(h)=g_+\}$, the matrix of $T$ is a block
diagonal one.

Similarly, the matrix of $(T^*)^{-1}$ is a block diagonal one with regard to the
decomposition $\coprod_{g_-\in G_-}\{h^*\in G^*:\beta_-(h)=g_-\}$. Therefore the
matrix of $T$ is also block diagonal with regard to the
decomposition $\coprod_{g_-\in G_-}\{h\in G:\beta_-(h)=g_-\}$. This is equivalent
to $T\{g\}$ is a linear combination of terms $\{h\}$ with $h_-=g_-$.

Thus we conclude that $T\{g\}$ is a linear combination of terms $\{h\}$ with
$h_+=g_+$ and $h_-=g_-$. There is only one such term, namely $\{g\}$. Therefore
$T\{g\}$ is a scalar multiple of $\{g\}$.

\hfill$\Box$

The fact that $g\mapsto \tau(g)g$ is a Hopf algebra isomorphism from $H_{\epsilon}$
to $H_0$ means exactly the following equalities
\begin{eqnarray}
\lambda(g,h)^{\epsilon} & = & 
\frac{\tau(g)\tau(h)}{\tau(gh_-)}
\hspace{5mm}\mbox{in case $\bar{g}_+=h_+$} \nonumber \\
\mu(g,h)^{\epsilon} & = & 
\frac{\tau(g_+h)}{\tau(g)\tau(h)}
\hspace{5mm}\mbox{in case $g_-=\bar{h}_-$} \nonumber \\
\nu(g)^{\epsilon} & = & \frac{\tau(g)}{\tau(g^{-1})} \nonumber
\end{eqnarray}
The scaling function in Lemma \ref{lem55} provides us with a function $\tau$
satisfying such equalities. Since all the structure constants $\lambda$, $\mu$,
$\nu$ are positive, this implies that the positive scaling $g\mapsto
|\tau(g)|^{\frac{1}{\epsilon}}$ satisfies the similar equalities without
$\epsilon$ power. In other words,
\[
g\mapsto |\tau(g)|^{\frac{1}{\epsilon}}g:\hspace{5mm} H_1\ra H_0
\]
is an isomorphism of Hopf algebras.

\section{Hopf Algebra in Correspondence Category}

In this section, we study the set-theoretical version of Theorem 1.
We will assume all the sets are finite.

As pointed out in the introduction, if we restrict ourselves to the usual category of
finite sets, with maps as morphisms, then we will end up with only finite groups.
Therefore we enlarge the category to include correspondences, so that we may get a more
interesting class of Hopf algebras.

A {\em correspondence} from a set $X$ to a set $Y$ is a subset $F$ of $X\times Y$.
Given a correspondence $F$ from $X$ to $Y$ and a correspondence from $Y$
to $Z$, we have the composition
\[
F\circ G=\{(x,z): \mbox{there is $y\in Y$, such that $(x,y)\in F$ and 
$(y,z)\in G$}\}.
\]
With correspondences as morphisms, the finite sets form a category, which we call the
{\em correspondence category}. With the usual product of sets as tensor product, a one
point set as an identity object, and the diagonal $I_X=\{(x,x): x \in X\}\sub X\times
X$ as the identity morphism, the correspondences category becomes a
{\em monoidal category} as defined in \cite{Mac}. 

One may introduce the concept of algebra (called monoid in \cite{Mac}) in any
monoidal category. By reversing the directions of arrows in all diagrams used in
defining algebra, we have the concept of coalgebra. In the correspondence category,
the product of two algebras is also an algebra. A {\em correspondence bialgebra}
consists of an algebra structure and a coalgebra structure, such that the coalgebra
map is an algebra morphism. A {\em correspondence Hopf algebra} is a correspondence
bialgebra with an antipode making the usual diagrams commutative.

The purpose of this section is to classify correspondence Hopf algebras.

It is not convenient to work with correspondences directly. So we take the
following Boolean viewpoint. First we observe that a correspondence 
$F\sub X\times Y$ from $X$ to $Y$ gives rise to a map on the power sets
\[
P_F: P(X)\ra P(Y), 
\]
\[
P_F(A)=
\{b\in Y: \mbox{There is $a\in A$, such that $(a,b)\in F$}\}.
\]
The map $P_F$ preserves union and sends empty set to empty set. Conversely, if a map
$P: P(X)\ra P(Y)$ preserves union and sends empty set to empty set, then we have
$P=P_F$ for the unique $F$ given by $F=\{(a,b): b\in P(\{a\})\}$.

Power sets can be considered as modules over the simplest {\em Boolean algebra}
${\bf B}$. ${\bf B}$ consists of two elements $0$ and $1$, and has the
commutative sum and product
\[
0+0=0,\hspace{5mm}
1+0=1+1=1,\hspace{5mm}
0\cdot 0=1\cdot 0=0,\hspace{5mm}
1\cdot 1=1.
\]
${\bf B}$ satisfies all the axioms of a ring except the existence of additive inverse.
A ${\bf B}$-{\em module} is a set $M$ with addition $+$, a special element $0$, and a
scalar multiplication ${\bf B}\times M\ra M$, such that 
\begin{enumerate}
\item $+$ is commutative, associative, and $0+x=x$ for all $x\in M$;
\item $1\cdot x=x$, $0\cdot x=0$;
\item $(b_1+b_2)\cdot x=b_1\cdot x+b_2\cdot x$, \hspace{5mm}
$b\cdot(x_1+x_2)=b\cdot x_1+b\cdot x_2$.
\end{enumerate}
It is straightforward to check that, under the first two conditions, the third
condition means exactly
\begin{equation}\label{bool}
x+x=x,\hspace{5mm}\mbox{for all $x\in M$}.
\end{equation}
In fact, a ${\bf B}$-module is equivalent to a lattice.
In view of our assumption that all sets are finite, we will also assume all 
${\bf B}$-modules are finite.

With the obvious notion of homomorphisms, (finite) ${\bf B}$-modules form a category. 
The power set $P(X)$ is a ${\bf B}$-module by
\[
A+B=A\cup B,\hspace{5mm}
0=\emptyset.
\]
It is easy to see that 
\[
P: X\mapsto P(X), F\mapsto P_F
\]
is a functor from the correspondence category to the ${\bf B}$-module category.

A subset $B$ of a ${\bf B}$-module $M$ is called a {\em basis} if every $x\in M$ is a
unique linear combination of elements in $B$. A ${\bf B}$-module is called {\em free}
if it has a basis. Clearly, the collection $\{\{x\}:x\in X\}$ is a basis of $P(X)$.
The following lemma (especially the third property) implies $P$ is a one-to-one
correspondence between finite sets and finite free ${\bf B}$-modules.

\begin{lem}\label{boolfree}
Let $B$ be a basis of a ${\bf B}$-module $M$. Let $b\in B$. 
\begin{enumerate}
\item If the $b$-coordinate of $x\in M$ is $1$, then for any $y\in M$, the
$b$-coordinate of $x+y$ is also $1$;
\item If $b=x_1+\cdots+x_k$ for some $x_i\in M$, then each $x_i=b$ or $0$;
\item The basis of $M$ is unique.
\end{enumerate}
\end{lem}

\noindent{\em Proof}: In the first property, we have $x=b+\xi$, where $\xi$
is a sum of elements in $B-\{b\}$. We also have two possibilities for $y$:
\[
y=b+\eta, \quad\mbox{or}\quad y=\eta,
\]
where $\eta$ is a sum of elements in $B-\{b\}$. Then by (\ref{bool}), we have
\[
x+y=\left\{\begin{array}{ll}
b+b+\xi+\eta & \mbox{if $y=b+\eta$} \\
b+\xi+\eta & \mbox{if $y=\eta$}
\end{array}\right.
=b+\xi+\eta,
\] 
and $\xi+\eta$ is still a sum of elements in $B-\{b\}$. This proves the first property.

In the second property, if $x_1$ is neither $b$ nor $0$, then we have $b'\in B-\{b\}$,
such that the $b'$-coordinate of $x_1$ is $1$. By taking $x=x_1$, $y=x_2+\cdots+x_k$,
and $b=b'$ in the first property, we see that the $b'$-coordinate of $b$ is also $1$.
The contradiction shows that $x_1=b$ or $0$. The proof for $x_i$ is similar.

If we take the expression $b=x_1+\cdots+x_k$ in the second property as the expansion
of $b$ (in one basis $B$) in terms of another basis (consisting of $x_1,\cdots,x_k$
and some others), then we immediately see that $B$ is the only basis of $M$.

\hfill$\Box$

The tensor product of ${\bf B}$-modules can be defined in the same way as modules over
commutative rings. With ${\bf B}$ as the unit object, ${\bf B}$-modules then form
a monoidal category. Since we clearly have $P(X)\otimes P(Y)=P(X\times Y)$ and
$P(\mbox{one point})={\bf B}$, $P$ is a functor of monoidal categories.

The concepts such as ${\bf B}$-algebra, ${\bf B}$-coalgebra,
and ${\bf B}$-Hopf algebra can be introduced just like the corresponding
concepts over commutative rings. Since $P$ is a functor of monoidal categories, $P$
sends correspondence Hopf algebras to free ${\bf B}$-Hopf algebras. Moreover,
because of the third property in Lemma \ref{boolfree}, $P$ is a one-to-one
correspondence between correspondence Hopf algebras and free ${\bf B}$-Hopf algebras.

Motivated by the classification of positively based Hopf algebras over ${\bf C}$, we
may introduce free ${\bf B}$-Hopf algebras $H(G;G_+,G_-)$ as in Section 2 and
conjecture that these are all the free ${\bf B}$-Hopf algebras. This is indeed the
case. 

\begin{th}\label{last}
Any free ${\bf B}$-Hopf algebra is isomorphic to $H(G;G_+,G_-)$ for a unique group
$G$ and a unique factorization $G=G_+G_-$.
\end{th}

\noindent{\em Proof}: The theorem may be proved by basically repeating the contents of
Section 2. There is no need to do rescaling (and all the positive coefficients in
Section 2 can be replaced by 1) since the only nonzero element of ${\bf B}$ is $1$.

Instead of repeating all the details, we discuss some key points involved in Section 2,
and prove only one lemma to illustrate that the whole section works in the category of
${\bf B}$-modules.

There are two key technical points in the proofs in Section 2:  One is the use of
duality, opposite, and co-opposite in translating a proven lemma into another lemma in
similar or dual settings. The other is the use of positivity.

The dual of ${\bf B}$-modules can be defined in the same way
as modules over commutative rings, and satisfies similar properties such as
$(M^*)^*=M$. The opposite and co-opposite of ${\bf B}$-Hopf algebras can be introduced
similarly, as far as the antipode is invertible. 

As for the positivity, the role can be replaced by the three properties in Lemma
\ref{bool}. Here we only illustrates how this works by proving Lemma
\ref{lem1}. Assume $B$ is the basis of a ${\bf B}$-Hopf algebra $H$. The identity
element $1$ of $H$ is a summation $1=d_1+\cdot+d_k$ of basis elements. Then we have 
$d_i=d_i\cdots 1=d_id_1+\cdots+d_id_k$. Because $d_i$ is a basis element, by the
second property in Lemma \ref{bool} we see that $d_id_j=0$ or $d_i$. Similarly,
from $d_j=1\cdot d_j$, we see that $d_id_j=0$ or $d_j$. Then we may conclude
$d_id_j=0$ when $i\neq j$ and $d_id_i=d_i$, as before. 

Similarly, from $\Delta (d_1)+\cdots+\Delta (d_k)=\Delta (1)=1\otimes 1=
\sum_{i,j}d_i\otimes d_j$, we may prove ${\bf B}d_1+\cdots+{\bf B}d_k$ is
closed under $\Delta$. From $S(d_1)+\cdots+S(d_k)=S(1)=1=d_1+\cdots+d_k$, we may
prove ${\bf B}d_1+\cdots+{\bf B}d_k$ is closed under $S$. Therefore 
${\bf B}d_1+\cdots+{\bf B}d_k$ is a Hopf subalgebra of $H$. The standard
argument for classifying commutative Hopf algebras also applies to ${\bf B}$-modules.
Then we know ${\bf B}d_1+\cdots+{\bf B}d_k$ must be the dual
${\bf B}$-Hopf algebra of the group algebra ${\bf B}G_+$. This completes the proof of
the Lemma \ref{lem1} for ${\bf B}$-Hopf algebras.

\hfill$\Box$

By the one-to-one correspondence between correspondence Hopf algebras and free ${\bf
B}$-Hopf algebras, we may translate Theorem \ref{last} into the following
classification of correspondence Hopf algebras.

\begin{th}\label{correspondence}
Any finite correspondence Hopf Algebra is a finite group $G$, with a unique
factorization $G=G_+G_-$, such that 
\begin{eqnarray}
m & = & 
\{(g_-h_+,h_+k_-,g_-h_+k_-): h_+\in G_+, g_-,k_-\in G_-\}
\sub (G\times G)\times G   \nonumber \\
\Delta  & = & 
\{(g_+h_-k_+,g_+h_-,h_-k_+): g_+,k_+\in G_+, h_-\in G_-\}
\sub G\times (G\times G)   \nonumber \\
1 & = & \{(pt,g_+): g_+\in G_+\}\sub \{pt\}\times G \nonumber \\
\epsilon & = & \{(g_,pt): g_-\in G_-\}\sub G\times\{pt\} \nonumber \\
S & = & \{(g,g^{-1}): g\in G\}\sub G\times G  \nonumber
\end{eqnarray}
\end{th}

\hfill$\Box$

\end{document}